\documentclass[11pt]{article}

\usepackage{times}
\usepackage{bm}
\usepackage{graphicx}
\usepackage{amsfonts}
\usepackage{amsmath}
\usepackage{amssymb}
\usepackage{bm}
\usepackage{enumerate}
\usepackage{booktabs,fixltx2e,subfigure} 
\usepackage[flushleft]{threeparttable}
\usepackage{color}
\usepackage[round]{natbib}
\usepackage{setspace,cite}

\newcommand{\disp}{\displaystyle}

\allowdisplaybreaks

\begin{document}
\title{
\vspace{-0.75in}
\textbf{\LARGE \sc
On Construction of Higher Order Kernels Using Fourier Transforms and Covariance Functions
\vspace{0.5in}
}}\small
 \small\author{
 \hspace{-0.2in}   \parbox[t]{0.425\textwidth}{{\small \sc  Soumya Das}
  \\ {\footnotesize CEMSE Division (Statistics)\\ King Abdullah University of Science and Technology \\  KSA \\ sundew122436@gmail.com}}
    \hspace{.0005\textwidth}
    \parbox[t]{0.35\textwidth}{{\small \sc  Subhajit Dutta}
 \thanks{Research has been partially supported by ECR/2017/000374, Science \& Engineering Research Board (SERB), Department of Science and Technology, Government of India.}
  \\ {\footnotesize Department of Mathematics and Statistics\\ Indian Institute of Technology Kanpur\\ India \\  tijahbus@gmail.com}}
    \hspace{.0005\textwidth}
  \parbox[t]{0.4\textwidth}{{\small \sc Radhenduhska Srivastava }
 \thanks{Research supported by INSPIRE fellowship, Department of Science and Technology, Government of India and a seed grant from IIT Bombay.}
  \\ {\footnotesize Department of Mathematics\\ Indian Institute of Technology Bombay\\  India \\ rsrivastava22@gmail.com}}
\vspace{0.25in}
}

\normalsize

\maketitle
\begin{abstract}
\small
In this paper, we show that a suitably chosen covariance function of a continuous time, second order stationary stochastic process can be viewed as a symmetric higher order kernel. This leads to the construction of a higher order kernel by choosing an appropriate covariance function. 
An optimal choice of the constructed higher order kernel that partially minimizes the mean integrated square error of the kernel density estimator is also discussed.
\end{abstract}

\vspace{0.25in}
and function.
\normalsize

\newpage

\section{Introduction} \label{introduction}

Given univariate independent and identically distributed (i.i.d.) random variables $X_1,\ldots,X_n$ with a common probability density function (pdf) $f(x)$ at a fixed point $x \in \mathbb{R}$, a commonly used nonparametric estimator of $f(x)$ is the kernel density estimator given by (see, e.g., \citet{Silverman_1986})
\begin{equation}\label{k_est}
f_n(x,h)=\frac 1 n \sum_{i=1}^n K_h(x - X_i),
\end{equation}
where $K_h(u)=h^{-1} K(u/h)$, $K$ is an even real-valued function referred as the kernel function, and $h>0$ is the smoothing parameter. It is well known that under fairly general conditions on $f$, the kernel density estimator $f_n$ is a consistent estimator of $f$ when $h\to 0$ together with $nh\to \infty$ as $n\to \infty$ (see, e.g., \citet{Silverman_1986}). Moreover, the rate of convergence of the bias of $f_n$ to zero can be sharpened if the underlying density function $f$ is smooth. In particular, if $f$ is  $(p+1)$ times continuously differentiable with bounded $f^{(p+1)}$, then by using Taylor's series expansion, we obtain
\begin{equation}\label{bias_fn}
E[f_n(x,h)] = \int_{-\infty}^{\infty} K(u)f(x-hu)du = f(x) + \sum_{l=1}^p \frac{(-h)^{l}f^{(l)}(x)}{l!} \int_{-\infty}^{\infty} u^{l}K(u)du + o(h^{p}),
\end{equation}
where $f^{(l)}$ denotes the $l^{\scriptsize \mbox{th}}$ derivative of $f$ for $l=1,\ldots,p$. Further, the fastest possible rate of the bias of $f_n$ could be $O(h^p)$, if the chosen kernel function $K$ satisfies
$\disp \int_{-\infty}^{\infty} u^{j} K(u)du=0$ for $j=1,\ldots,(p-1)$.
We now discuss a closely related notion. 
Throughout the paper, $L^{p}(\mathbb{R})$ denotes the class of functions $g$ with $\disp \int_{-\infty}^{\infty} |g(u)|^p du < \infty$ for $p \geq 1$.
\vskip 5pt

{\bf Definition 1.} A real-valued function $K\in L^{1}(\mathbb{R})$ is said to be a kernel of order $p$, if it satisfies the following conditions.
\begin{enumerate}[(a)]
\item $\disp \int_{-\infty}^{\infty} K(u) du=1$,
\item $\disp \int_{-\infty}^{\infty} u^j K(u)du=0$ for $j=1,2,\ldots,(p-1)$ and $\disp \int_{-\infty}^{\infty} u^p K(u)du\ne0$.
\end{enumerate}


\vspace{0.05in}
Several researchers have explored the construction of higher order kernels to improve the rate of convergence of the density estimator (see \citet{Parzen_1962, Rosenblatt_1971}). \citet{WS_1990} considered an approach composed of polynomials and a second order kernel function to construct higher order kernels. The twicing technique of \citet{Abdous_1995} when applied to a kernel $K$ of order $p$ leads to the kernel function $2K - K * K$ 
of order $2p$, 
where $g_1*g_2(t)= \disp \int_{-\infty}^{\infty} g_1(u) g_2(t-u) du $ is the convolution of the functions $g_1, g_2 \in L^1(\mathbb{R})$.
Generalizations for higher order kernels based on twicing approaches are usually quite complicated, and may not yield explicit expressions. Construction of such higher order kernels involve an iterative twicing procedure as well as generalized jackknifing.
\citet{GHU_2003} provides a good review of the existing procedures of constructing higher order kernels (also see \citet{Scott_2015}).
\citet{HM_1988} and \citet{UU_2012} studied higher order kernel based on Fourier transforms.

The asymptotic order of the bias of $f_n$ is limited by the order of the kernel when the underlying pdf $f$ is infinitely differentiable.
\citet{Davis_1975} studied the behavior of $f_n$ by using the {\it sinc kernel} defined as $K(u) = (\pi u)^{-1}\sin(u)$ which is a kernel of `infinite' order. 
\citet{Devroye_1992} constructed super-kernels which are `infinite' order kernels (see also \citet{PR_1999}). 
If the underlying density function $f$ is $p$ times continuously differentiable with bounded $p^{\scriptsize \mbox{th}}$ derivative, then the rate of the bias of a super-kernel density estimator $f_n$ is $o(h^p)$.
\citet{GHU_2002} concluded that the {\it sinc kernel} density estimator is preferable compared to super-kernel density estimators in several situations.
Recently, \citet{Chacon_2014} have used the {\it sinc kernel} for smooth distribution function estimation.

In this article, we aim to construct an explicit family of higher order kernel for a given even order. In Section \ref{Construct_higher order kernel}, we provide a general methodology to construct higher order kernel using Fourier transforms. This method requires us to compute the Fourier transformation of an appropriately chosen function, which may not be easy to obtain explicitly. We also view this construction as a covariance function of a continuous time, second order stationary stochastic process. Using this view point, we construct a family of kernels motivated from Shannon's formula \citep{Shannon_1949} based on the {\it sinc function} (defined in Section \ref{higher order kernelII}).
In Section~\ref{higher order kernel_MISE}, we obtain an expression for the MISE of the kernel density estimator for the chosen kernel and provide a guideline to choose the kernel from the constructed family of kernels.
In Section~\ref{Num_higher order kernel}, we compare the mean integrated square error (MISE) of the kernel density estimator obtained from our chosen higher order kernel with some regular kernels and the {\it sinc kernel} using some numerical examples.
All theoretical proofs are given in the Appendix, and additional results are stated in a Supplementary file.

\section{Construction of a Higher Order Kernel Using Fourier Transforms} \label{Construct_higher order kernel}

The following theorem gives a general method to construct symmetric higher order kernels.

\medskip\noindent
{\bf Theorem 1.}
{\it Let $G$ be even, $q~(\ge 1)$ times continuously differentiable and compactly supported function such that $G(0)=1$ and $G^{(j)}(0)=0$ for $j=1,2,\ldots,(2q-1)$ and $G^{(2q)}(0)\ne0$.
Let $K$ be the Fourier transform of $G$ defined as
\begin{eqnarray*}
K(u)&=&\int_{-\infty}^{\infty} G(t) e^{-2\pi i u t} dt,
\end{eqnarray*}
where $i=\sqrt{-1}$. Then, $K$ is a symmetric kernel of order $2q$.}\\
\citet{HM_1988} also considered kernels constructed using Fourier transform of functions, but they considered a specific choice of $G$. The proof of Theorem 1 can be established easily. For the sake of completeness, we have provided the proof in the Supplementary material.

If the function $G$ is positive, then the constructed higher order kernel $K$ can be viewed as covariance function of a continuous time second order stationary stochastic process with power spectral density $G$. For example, the covariance functions corresponding to the following $2q$ degree polynomial power spectral densities, are the higher order kernels of order $2q$. Consider the following examples:
\begin{enumerate}
\item $G_1(t)=(1 -4t^2)^{q} I_{[-1/2,1/2]}(t)$,
\item $G_2(t)=(1-(2t)^{2q}) I_{[-1/2,1/2]}(t)$,
\end{enumerate}
\noindent
where $I_A$ is the indicator function associated with the set $A$. The higher order kernels (equivalently, the covariance functions) corresponding to the spectral density $G_1$ is listed in Table~1 for different orders.
\begin{table}[ht]
\caption{Higher order kernels (covariance function) based on the spectral density $G_1$} 
\begin{center}
{\scriptsize
\vspace{0.1in}
\begin{tabular}{ccc}
\hline \hline
$p$ & $K(0)$ & $K(u)^*$
 \\ [1ex]
\hline
2 &2/3 & $\{16\sin(u/2) - 8u\cos(u/2)\}/{u^3}$ \\ [1ex]
4 &8/15 & $\{-384u \cos(u/2) + 768 \sin(u/2) - 64u^2 \sin(u/2)\}/{u^5}$ \\ [1ex]
6 &16/35 & $\{92160 \sin(u/2) - 46080 u \cos(u/2) + 768 u^3 \cos(u/2) - 9216 u^2 \sin(u/2)\}/{u^7}$\\[1ex]
8 &128/315 & $\{20643840 \sin(u/2) - 10321920 u \cos(u/2) + 245760 u^3 \cos(u/2) - 2211840 u^2 \sin(u/2) + 12288 u^4 \sin(u/2)\}/{u^9}$\\
[1ex]
\hline
\end{tabular}
}
\end{center}
\label{t1}
\footnotesize{* The function $K$ is obtained by using the function {\tt fourier} in MATLAB.}
\end{table}


\subsection{Higher Order Kernel Using Truncated Covariance Function} \label{higher order kernelII}

Given a general $G$ satisfying conditions of Theorem 1, the expression of $K$ may not be easily accessible (e.g., the spectral density $G_2$). In this section, we construct a higher order kernel which has a closed form expression.
Suppose that the support of spectral density $G$ is $[-1/2,1/2]$, then using Shannon's formula \citep{Shannon_1949}, the covariance function $K$ can be expressed as
\begin{eqnarray}\label{kernel_cons}
K(u)&=&\sum_{j=-\infty}^\infty K\left(j\right)\mbox{sinc}\left(\pi \left(u-j\right)\right),
\end{eqnarray}
where $\mbox{sinc}(u)={\sin u}/u$ if $u \neq 0$, and $1$ if $u=0$.
The representation \eqref{kernel_cons} shows that the function $K(u)$ can be reconstructed from the sequence $\{K(j), \ j=\ldots,-1,0,1,\ldots\}$. Further, if
$\sum_{j=-\infty}^{\infty}\left| K(j) \right|<\infty$, then
\begin{eqnarray}\label{spec_den}
G(t) &=&\int_{-\infty}^{\infty} K(u)e^{2\pi i t u} du=\sum_{j=-\infty}^{\infty}K\left(j\right) e^{-2\pi i j t }\quad 1_{[-\frac12,\frac12]}(t).
\end{eqnarray}
The equality in equation (\ref{kernel_cons}) holds in $L^2(\mathbb{R})$ sense, i.e., $\disp \int_{-\infty}^{\infty} (K_N(u)-K(u))^2 du \to 0$ as $N \to \infty$, where $K_N(u) = \disp \sum_{j=-N}^N K(j) \; \mbox{sinc}(\pi(u-j))$ (see \citet{Shannon_1949}).
Additionally, if $\disp \sum_{j=1}^{\infty} \frac{|K(j)|}{j}<\infty$, then the equation in \eqref{kernel_cons} holds pointwise as well (see Lemma 1 in the Supplementary material).

In view of equations \eqref{kernel_cons} and \eqref{spec_den}, the condition $G(0)=1$ is equivalent to $\sum_{j=-\infty}^\infty K(j)=1$. Similarly, conditions $G^{(2r)}(0)=0$ corresponds to $\sum_{j=-\infty}^\infty j^{2r}K(j)=0$ for $r=1,2,\ldots,q-1$ and $G^{(2q)}(0) \neq 0$ corresponds to $\sum_{j=-\infty}^\infty j^{2q}K(j) \neq 0$.
If we truncate the series in equation \eqref{kernel_cons} by choosing $K(j)=0$ for all $|j|>q$, then the conditions on $G$ stated in Theorem~1 in terms of sequence $\{K(j), j=\ldots,-1,0,1,\ldots\}$ reduces to the following
\begin{eqnarray}
            &&\sum_{j=-q}^q K(j)=1,\qquad  \sum_{j=-q}^q j^{2r}K(j)=0 \mbox{ for } r=1,\ldots,(q-1),\label{system}\\
 \mbox{and} && \sum_{j=-q}^q  j^{2q}K(j) \neq 0.\label{non_zero}
\end{eqnarray}
A solution to the system of linear equation \eqref{system} that satisfies \eqref{non_zero} leads to a higher order kernel of order $2q$ by using equation \eqref{kernel_cons}.
Let $K(0)=\alpha$ and by using the symmetry of $K$, i.e., $K(-u)=K(u)$, the system of linear equation \eqref{system} reduces to the following:
\begin{equation} \label{solve}
\sum_{j=1}^q K(j)=(1-\alpha)/2, \ \  \sum_{j=1}^q j^{2} K(j)=0, \ \ldots, \ \ \sum_{j=1}^q j^{2(q-1)} K(j)=0.
\end{equation}

\noindent
{\bf Theorem 2.} {\it
If $\alpha \neq 1$, the system of linear equations stated in \eqref{solve} has the unique solution
\begin{eqnarray*}
&&K(0)=\alpha \mbox{ and } K(j)=\, \frac{1-\alpha}{2} \left(\frac{q!}{j}\right)^2 \frac{1}{\prod \limits_{\substack{l=1 \\ l \neq j}}^q (l^2-j^2)} , \hspace{0.5cm}  \forall\,j=1,\ldots,q,\\
\mbox{Further,}&&  \sum_{j=-q}^q  \  \ j^{2q}K(j)= (1-\alpha) (q!)^{2} (-1)^{q+1}   \neq 0.
\end{eqnarray*}
}
\vskip-20pt
\noindent
In view of Theorem~2, the kernel function $K$ obtained by using equation \eqref{kernel_cons}, i.e.,
\begin{eqnarray}\label{kernel_final}
K(u)&=& \disp \sum_{j=-q}^q K\left( j \right)\mbox{sinc}\left(\pi \left(u-j\right)\right),
\end{eqnarray}
where $\{K(j), j=0,\pm 1,\ldots,\pm q\}$ is as stated in Theorem~2, is a higher order symmetric kernel of order $2q$.

\section{Mean Integrated Square Error of Density Estimator} \label{higher order kernel_MISE}

The expression for the kernel density estimator is $f_n(x,h)=\frac 1 n \sum_{i=1}^n K_h(x - X_i)$,
where $K_h(u)=h^{-1}K(u/h)$ and $h>0$. The mean integrated square error (MISE) of $f_n(\cdot,h)$ is defined as
$$ MISE(f_n(\cdot,h))= E_f \int_{-\infty}^{\infty} \bigl \{f_n(x,h) - f(x) \bigr \}^2 dx.$$
Since the kernel function $K$ is the Fourier transform of the function $G$ (see Theorem 1), we obtain a simplified expression of the MISE as follows:
\begin{equation}\label{MISE}
MISE(f_n(\cdot,h)) = \frac1{nh}\int_{-\infty}^{\infty} G^2(t)dt + \int_{-\infty}^{\infty} \left[G(ht)-1\right]^2 |\phi_f(t)|^2 dt - \frac 1 n \int_{-\infty}^{\infty} G^2(ht)|\phi_f(t)|^2 dt,
\end{equation}
where $\phi_f(t)=\disp \int_{-\infty}^{\infty} f(u) e^{2\pi i tu} du$ and $G(t) =\disp \int_{-\infty}^\infty K(u)e^{2\pi i t u } du$ (see \citet{UU_2012}).

\vspace{0.1in}
\noindent
{\bf Theorem 3.} {\it
Suppose that the density function $f$ is $2p$ times differentiable, where $p$ is an even number. Then, the MISE of the density estimator corresponding to the kernel function $K$ of order $p$, as constructed in \eqref{kernel_final}, is given by
\begin{eqnarray*}
\lim_{n\to\infty}n^{\frac{2p}{2p+1}}MISE(f_n(\cdot,h))&=& \int_{-\infty}^{\infty} G^2(t)dt + \frac{(G^{(p)}(0))^2}{(p!)^2} \int_{-\infty}^{\infty} t^{2p} |\phi_{f}(t)|^2 dt.
\end{eqnarray*}

\noindent
If we minimize the first term of the MISE expression stated above over $\alpha$, then the choice of $\alpha$ is $\disp \frac{C_p}{1+C_p}$, where $C_p= \disp \frac 1 2 \sum_{l=1}^p \biggl ( \frac{(p!)^2}{l^2 \prod \limits_{j=1, j \neq l}^p (j^2-l^2)} \biggr )^2$.
}
\vspace{0.1in}

Although the kernel function $K$ obtained by choosing the function $G$ as per Theorem 1 is a higher order kernel of order $p$, the explicit form of this kernel as in Theorem 2 holds in $L^2(\mathbb{R})$ sense and pointwise 
but not in $L^1(\mathbb{R})$ sense as the $\it{sinc~function}$ is not integrable. Thus, the kernel function constructed in Theorem 2 is not an integrable function but square integrable, and the resulting density estimator is not a valid density estimator. Further, the density estimator corresponding to the higher order kernel constructed in Theorem 2 is neither non-negative nor integrable.

To rectify this problem, we use the proposal of \citet{GHU_2003} and define
\begin{equation*}
{\tilde f}_n(x,h)=\max\{0,f_n(x,h)-\xi\},
\end{equation*}
where the constant $\xi$ is chosen such that $\disp \int_{-\infty}^{\infty} {\tilde f}_n(x,h) dx=1$.
This correction ensures that ${\tilde f}_n(x,h)$ is non-negative, integrates to one and it also follows that the MISE of this new version (say, $MISE({\tilde f}_n(\cdot,h))$) is lower than $MISE({f}_n(\cdot,h))$ (see Theorem 1 in \citet{GHU_2003}).

\section{Finite Sample Performance}\label{Num_higher order kernel}

In this section, we compare the MISEs of the conventional (with the Gaussian kernel) density estimator, the usual sinc density estimator, and the two density estimators proposed in Section \ref{Construct_higher order kernel} 
for finite values of the sample size (say, $n$). We have considered three values for the sample size $n = 50$ (small), $n = 250$ (medium) and $n = 500$ (large) over $100$ simulated samples. The samples are drawn from the $N(0,0.1)$ distribution, the gamma distribution with shape and rate parameters both equal to $2$,
the $l_p$-symmetric (i.e., $f(x)= (p/2\Gamma{(p)})\exp(-|x|^p)$) with $p=3$ distribution, and the Fej{\'e}r-de la Vall{\'e}e Poussin (FVP) density (i.e., $f(x)=(2/\pi)(x^{-1}\sin(x/2))^2$).
For numerical experiments, we have used the corrected versions of the sinc and the Fourier based density estimators. The bandwidth selection approaches for the competing methods were different. We have used the function {\tt bw.nrd} in the R package {\tt stats} for the Gaussian kernel, and  ${(\log(n + 1))}^{-1/2}$ for the sinc kernel \citep{GHU_2002}. For the proposed methods, we have taken $n^{-1/(2p + 1)}$ when the kernel is of order $p$. The grid over which an estimator evaluated is $1001$ equi-spaced points in the interval $[-5,5]$. Results of the average MISEs are reported in Table \ref{Table.Num}. The minimum MISE is marked in {\bf bold}, while the second best is marked in {\it italics}.

It is clear from Table \ref{Table.Num} that the density estimator with the Gaussian kernel performs quite well for the first two examples, while the estimator based on $G_1$ (a higher order kernel of order $2$) stated in Table \ref{t1} yields the best performance in the next two examples.
In the third and fourth examples, the $l_3$-symmetric density is not differentiable at the point $0$, while the FVP density is quite wriggly, respectively.
The overall performance of the kernel based on the truncated {\it sinc function} (say, tsinc) of order $2$ (see equation \eqref{kernel_final} of Section \ref{higher order kernelII}) is quite competitive.
 
\begin{table}[h]
	\caption{Estimated MISE (with standard error in brackets) for varying sample sizes} 
	\vspace{0.2in}
	\centering 
	\begin{tabular}{l c cccc} 
		\hline\hline 
		Distribution & $n$ & Gaussian & sinc & $G_1$ & tsinc
		\\ [0.5ex]
		\hline 
		&$50$ & {\bf 0.0032 (0.0002)} & 0.0290 (0.0004) & 0.0386 (0.0001) & {\it 0.0170 (0.0008)} \\[-1ex]
		\raisebox{1.5ex}{$N(0,0.1)$} & $250$
		& {\bf 0.0010 (0.0006)} &0.0215 (0.0002) & 0.0238 (0.0001) & {\it 0.0067 (0.0007)} \\[1ex] & $500$
		& {\bf 0.0006 (0.0003)} & 0.0191 (0.0002) & 0.0174 (0.0001) & {\it 0.0043 (0.0005)} \\[1ex] \hline
		&$50$ & {\bf 0.0034 (0.0002)} & 0.0081 (0.0007) & 0.0105 (0.0005) & {\it 0.0056 (0.0001)} \\[-1ex]
		\raisebox{1.5ex}{$G(2,2)$} & $250$
		& {\bf 0.0015 (0.0005)} & 0.0058 (0.0003) & 0.0059 (0.0005) & {\it 0.0026 (0.0005)} \\[1ex] & $500$
		& {\bf 0.0010 (0.0003)} & 0.0051 (0.0002) & 0.0050 (0.0005) & {\it 0.0018 (0.0003)} \\[1ex] \hline
		&$50$ & {\it 0.0083 (0.0002)} & 0.0107 (0.0001) & {\bf 0.0073 (0.0002)} & 0.0096 (0.0002) \\[-1ex]
		\raisebox{1.5ex}{$l_3(0,1)$} & $250$
		& {\it 0.0063 (0.0008)} & 0.0075 (0.0005) & {\bf 0.0048 (0.0002)} & 0.0070 (0.0007) \\[1ex] & $500$
		& {\it 0.0060 (0.0005)} & 0.0069 (0.0003) & {\bf 0.0040 (0.0001)} & 0.0064 (0.0004) \\[1ex] \hline	
		&$50$ & 0.0069 (0.0002) & {\it 0.0066 (0.0003)} & {\bf 0.0022 (0.0002)} & {0.0068 (0.0002)} \\[-1ex]
		\raisebox{1.5ex}{FVP} & $250$
		& 0.0060 (0.0008) & {\it 0.0056 (0.0001)} & {\bf 0.0010 (0.0006)} & {0.0057 (0.0001)} \\[1ex] & $500$
		& 0.0058 (0.0006) & {\it 0.0054 (0.0008)} & {\bf 0.0009 (0.0005)} & {0.0055 (0.0007)} \\[1ex] \hline		
	\end{tabular}
	\label{Table.Num}
\end{table}

\section{Appendix} \label{higher order kernel_Append}

\noindent
{\bf Proof of Theorem 2.}
The system of linear equation given in (\ref{solve}) is expressed in the matrix form as follows.
\begin{equation}\label{matrix_form}
 \left[ \begin{array}{ccccc}
1 & 1 & 1 & \ldots & 1 \\
1^2 & 2^2 & 3^2 & \ldots & q^2 \\
1^4 & 2^4 & 3^4 & \ldots & q^4 \\
\vdots & \vdots & \vdots & \vdots & \vdots \\
1^{2(q-1)} & 2^{2(q-1)} & 3^{2(q-1)} & \ldots & q^{2(q-1)}
\end{array} \right]
 \left( \begin{array}{c}
K(1) \\
K(2) \\
K(3) \\
\vdots \\
K(q)
\end{array} \right)
=
 \left( \begin{array}{c}
(1-\alpha)/2 \\
0 \\
0 \\
\vdots \\
0
\end{array} \right).
\end{equation}
\noindent
The $(i,j)$th element of a $q\times q$ Vandermonde matrix is defined as $b_j^{i-1}$, where $b_j, j=1,2,\ldots,q$ are non-zero real numbers and the determinant of this matrix is $\prod \limits_{1 \leq i < j \leq q} (b_j-b_i)$ \citep{RB_2000}. The coefficient matrix (say, $A$) on the left hand side of  \eqref{matrix_form} is a $q \times q$ Vandermonde matrix with $b_j=j^2$ for $j=1,\ldots,q$. The determinant of coefficient matrix is non-zero, thus system of linear equations \eqref{solve} has a unique solution.

It is enough to compute the first column of $A^{-1}$ to get the solution vector of \eqref{matrix_form}.
Let $(i,j)^{\text{th}}$ element of $A^{-1}$ is denoted by $a^{ij}$ for $i,j=1,\ldots,q$.
Recall that $A^{-1}=\, \disp \frac{Adj A}{det A}$, where the adjugate of $A$ is the transpose of the co-factor matrix $C$ of $A$, i.e., $Adj(A)=C^T$.
Here, the $(i,j)^{\text {th}}$ entry of $C$ is $C_{ij}=\, (-1)^{i+j}M_{ij}$ with $M_{ij}$ being the $(i,j)^{\text{th}}$ minor of A for $(i,j) \in \{1,\ldots,q\}^2$. Now,
\begin{align*}
a^{k1}&=(-1)^{k-1}\frac{\begin{vmatrix}
1^2 & 2^2 & \dots & (k-1)^2 & (k+1)^2 & \cdots & q^2\\
1^4 & 2^4 & \dots & (k-1)^4 & (k+1)^4 & \cdots & q^4\\
\vdots & \vdots & \ddots & \vdots &\vdots & \ddots & \vdots \\
1^{2(q-1)} & 2^{2(q-1)} & \dots & (k-1)^{2(q-1)} & (k+1)^{2(q-1)} & \cdots & q^{2(q-1)}
\end{vmatrix}}{\begin{vmatrix}
1 & 1 & \dots & 1 \\
1^2 &2^2 & \dots &q^2 \\
\vdots & \vdots & \ddots & \vdots \\
1^{2(q-1)} &2^{2(q-1)} & \dots &q^{2(q-1)}
\end{vmatrix}}\\
&=(-1)^{k-1} {1^22^2 \cdots (k-1)^2 (k+1)^2 \cdots q^2}\,  \frac{\prod \limits_{\substack{1\leq i < j \leq q \\ i \neq k , j \neq k }}(j^2-i^2)}{\prod \limits_{1\leq i < j \leq q} (j^2-i^2)}\\
&=(-1)^{k-1} \left(\frac{q!}{k} \right)^2 \, \frac{1}{{\prod \limits_{j=k+1}^q (j^2-k^2)}{\prod \limits_{i=1}^{k-1} (k^2-i^2)}} \\
&=\left(\frac{q!}{k}\right)^2 \frac{1}{\prod \limits_{\substack{j=1 \\ j \neq k}}^q (j^2-k^2)} , \hspace{5cm}  \forall\,k=1,2, \ldots ,q .
\end{align*}
Thus, the solution vector is given by
\begin{eqnarray} \label{sol_eq}
K(0)=\alpha, \mbox{ and } K(k)=\, \frac{1-\alpha}{2} \left(\frac{q!}{k}\right)^2 \frac{1}{\prod \limits_{\substack{j=1 \\ j \neq k}}^q (j^2-k^2)} \hspace{0.5cm}  \forall\,k=1,2,\ldots,q.
\end{eqnarray}
Now, for the solution vector
$$	\sum \limits_{k=-q}^{q} k^{2q} K(k) = 2 \sum \limits_{k=1}^{q} k^{2q} K(k) = (1-\alpha) (q!)^2 \displaystyle{\sum \limits_{k=1}^{q} k^{2(q-1)}\frac{1}{\prod \limits_{\substack{j=1 \\ j \neq k}}^q (j^2-k^2)}} = (1-\alpha)(q!)^2 \displaystyle{\sum \limits_{k=1}^q \frac{1}{\prod \limits_{\substack{j=1 \\  j \neq k}}^q (\frac{j^2}{k^2}-1)}}.$$
If $\alpha \neq 1$, then by Lemma 2 (stated in Supplementary material), we have the non-trivial identity $$\displaystyle{\sum \limits_{k=1}^q \frac{1}{\prod \limits_{\substack{j=1 \\  j \neq k}}^q (\frac{j^2}{k^2}-1)}} =(-1)^{q+1}.$$
This completes the proof. \hfill $\Box$

\vspace{0.1in}

\noindent
{\bf Proof of Theorem 3.} In view of \eqref{MISE}, the MISE of the density estimator $f_n(\cdot,h)$ is given by
\begin{eqnarray*}
MISE(f_n(\cdot,h))&=& \frac1{nh}\int_{-\infty}^{\infty} G^2(t)dt + \int_{-\infty}^{\infty} \left[ G(ht)-1\right]^2 |\phi_f(t)|^2 dt - \frac1n \int_{-\infty}^{\infty} G^2(ht)|\phi_f(t)|^2 dt\\
&=& A_1+A_2+A_3 ~(\text{say}),
\end{eqnarray*}
where $\phi_f(t)= \disp \int_{-\infty}^{\infty} f(u) e^{2\pi i tu} du$.
Clearly, $A_1=O((nh)^{-1})$.

\noindent
By using Taylor series and the conditions on $G$ stated in Theorem 1, we have
\begin{eqnarray}
G(ht) &=& 1+\frac{1}{p!}G^{(p)}(\xi) h^p t^p,
\end{eqnarray}
where $\xi\in (0,th)$. Now, the second term $A_2$ reduces to
\begin{eqnarray*}
A_2 &=& \frac{h^{2p}}{(p!)^2}\int_{-\infty}^{\infty} (G^{(p)}(\xi))^2 t^{2p} |\phi_{f}(t)|^2 dt.
\end{eqnarray*}
By using the inverse Fourier transformation of $\phi_f$ and $2p$ times differentiability of $f$, we have
\begin{eqnarray*}
\frac{d^{2p}}{dx^{2p}} f(x)&=& \frac{d^{2p}}{dx^{2p}}\int_{-\infty}^{\infty} e^{-2\pi i t x} \phi_f(t) dt=\int_{-\infty}^{\infty} (-2\pi i)^{2p} t^{2p} e^{-2\pi i t x} \phi_f(t) dt,
\end{eqnarray*}
where the interchange of the differentiation with the integral follows by Lebesgue's DCT. Thus, $t^{2p} \phi_f(t)$ is integrable. Since $|\phi_f(t)|\le 1$ and by using DCT, the function $t^{2p} |\phi_f(t)|^2$ is also integrable.
Again by applying Lebesgue's DCT, we get
\begin{eqnarray*}
\lim_{h\to 0} \frac{A_2}{h^{2p}}  &=& \frac{(G^{(p)}(0))^2}{(p!)^2} \int_{-\infty}^{\infty} t^{2p} |\phi_{f}(t)|^2 dt.
\end{eqnarray*}
We now turn to the third term $A_3$.
\begin{eqnarray*}
  A_3 &=&\frac1n \int_{-\infty}^{\infty} |\phi_{f}(t)|^2 dt+ \frac1n \frac{h^{2p}}{(p!)^2}\int_{-\infty}^{\infty} (G^{(p)}(\xi))^2 t^{2p} |\phi_{f}(t)|^2 dt+
  \frac{2h^{p}}{np!}\int_{-\infty}^{\infty} G^{(p)}(\xi) t^{p} |\phi_{f}(t)|^2 dt\\
  &=& A_{31}+ A_{32}+A_{33} ~(say),
\end{eqnarray*}
Again, by using Lebesgue's DCT, we obtain
$$\lim_{h\to 0} \frac{n}{h^{2p}} A_{32}  = \frac{(G^{(p)}(0))^2}{(p!)^2}\int_{-\infty}^{\infty}  t^{2p} |\phi_{f}(t)|^2 dt \mbox{ and }
\lim_{h\to 0} \frac{n}{h^{p}} A_{33}  = \frac{G^{(p)}(0)}{p!}\int_{-\infty}^{\infty}  t^{p} |\phi_{f}(t)|^2 dt.$$
Thus, we have $$A_3=O(n^{-1}).$$
The rate of convergence of the term $A_3$ is faster than that of the terms $A_1$ and $A_2$. Now, by combining the rates of $A_1, A_2$ and $A_3$, we obtain
$$MISE(f_n(\cdot,h)) = O(({nh})^{-1})+ O(h^{2p}).$$
We equate both the rates to get the optimal rate of convergence of MISE and this leads to the choice of $h=n^{-\frac{1}{2p+1}}$. Subsequently, MISE is given by
\vspace{-0.1in}
\begin{eqnarray*}
\vspace{-0.1in}
\lim_{n\to\infty}n^{\frac{2p}{2p+1}}MISE(f_n(\cdot,h))&=& \int_{-\infty}^{\infty} G^2(t)dt +\frac{(G^{(p)}(0))^2}{(p!)^2}\int_{-\infty}^{\infty} t^{2p} |\phi_{f}(t)|^2 dt.
\end{eqnarray*}

The MISE of the density estimator depends on the kernel function through $\disp \int_{-\infty}^{\infty} G^2(t)dt$ and $(G^{(p)}(0))^2$. For the kernel constructed in \eqref{kernel_final}, by using Lemma 2, we have $(G^{(p)}(0))^2= ((p/2)!)^4 (1-\alpha)^2/4$. Ideally, one should choose a kernel which minimizes $n^{2p/(2p+1)}MISE(f_n(\cdot,h))$ expressed as
$$ \int_{-\infty}^{\infty} G^2(t)dt + \frac{((p/2)!)^4 (1-\alpha)^2}{4 (p!)^2}\int_{-\infty}^{\infty} t^{2p} |\phi_{f}(t)|^2 dt.$$
Clearly, the second term depends on the unknown density function $f$.
Thus, we choose the function $G$ that minimizes  the first term $\disp \int_{-\infty}^{\infty} G^2(t)dt$. By using Plancherel's theorem, we have $\disp \int_{-\infty}^{\infty} G^2(t)dt = \sum_{j=-\infty}^\infty K^2(j)$.
For the kernel function $K$ constructed in \eqref{kernel_final}, we now have
\begin{eqnarray*}
\int_{-\infty}^{\infty} G^2(t)dt &=& K^2(0)+2 \sum_{l=1}^p K^2(l) = \alpha^2+2 \frac{(1-\alpha)^2}{4} \frac 1 2 \sum_{l=1}^p \bigg ( \frac{(p!)^2}{l^2 \prod \limits_{{j=1, j \neq l}}^p (j^2-l^2)} \bigg )^2 \\
& = & \alpha^2 +(1-\alpha)^2 C_p, \mbox{ where }C_p= \disp \frac 1 2 \sum_{l=1}^p \bigg ( \frac{(p!)^2}{l^2 \prod \limits_{{j=1, j \neq l}}^p (j^2-l^2)} \bigg )^2.
\end{eqnarray*}
If we minimize this expression over $\alpha$, the optimal solution turns out to be $\disp \frac{C_p}{1+C_p}$, which also satisfies $\disp \int_{-\infty}^{\infty} G^2(t)dt = \frac{C_p}{1+C_p}$.
This completes the proof. \hfill $\Box$

\small
\bibliographystyle{abbrvnat}
\bibliography{reference}

\end{document}